\documentclass[letterpaper, 10 pt, conference]{ieeeconf}

\IEEEoverridecommandlockouts

\usepackage{graphics}
\usepackage{epsfig}
\usepackage{times}
\usepackage{amsmath}
\usepackage{amssymb}
\usepackage{multirow}
\usepackage{array}
\usepackage{amsfonts}
\usepackage{subcaption}

\usepackage{amsthm}
\usepackage{todonotes}
\usepackage{comment}
\usepackage{cite}
\usepackage{booktabs}

\theoremstyle{definition}
\newtheorem{assumption}{Assumption}

\theoremstyle{plain}

\theoremstyle{definition}

\title{\LARGE \bf
Impact of Connected and Automated Vehicles in a Corridor}

\author{A M Ishtiaque Mahbub, {\itshape{Student Member, IEEE}}, Andreas A. Malikopoulos, {\itshape{Senior Member, IEEE}}, and \\ Liuhui Zhao, {\itshape{Member, IEEE}}
\thanks{This research was supported by the ARPAE's NEXTCAR program under the award number DE-AR0000796.}%
\thanks{The authors are with the Department of Mechanical Engineering, University of Delaware, Newark, DE 19716 USA (email: \tt\small{mahbub@udel.edu}; \tt\small{andreas@udel.edu}; \tt\small{lhzhao@udel.edu}.)} }

\begin{document}
\maketitle
\thispagestyle{empty}
\pagestyle{empty}

\begin{abstract}
Several approaches have been proposed in the literature that allow connected and automated vehicles (CAVs) to coordinate in areas where there is a potential conflict, for example, in intersections, merging at roadways and roundabouts. In this paper, we consider the problem of coordinating CAVs in a corridor consisting of several conflict areas where collision may occur. We derive a solution that yields the optimal control input, in terms of fuel consumption, for each CAV to cross the corridor under the hard safety constraints. We validate the effectiveness of the solution through simulation, and we show that both fuel consumption and travel time can be improved significantly.
\end{abstract}



\section{Introduction} \label{sec:in}

Connected and automated vehicles (CAVs) provide the most intriguing opportunity for enabling users to better monitor transportation network conditions and make better operating decisions. 
Several research efforts have been reported in the literature proposing different approaches on coordinating CAVs at different transportation scenarios, e.g., intersections, roundabouts, merging roadways, speed reduction zones, with the intention to improve traffic flow. In 2004, Dresner and Stone \cite{Dresner2004} proposed the use of the reservation scheme to control a single intersection. Since then, similar approaches have been reported in the literature to achieve safe and efficient control of traffic through intersections, e.g., \cite{Dresner2008, DeLaFortelle2010, Huang2012}. 

The objective of improving traffic flow through coordinating vehicles has been the focus in several papers \cite{Yan2009,Zhu2015a}. Kim and Kumar \cite{Kim2014} proposed an approach based on model predictive control that allows each vehicle to optimize its movement locally in a distributed manner with respect to any objective of interest. Most recently, \cite{Zhou2018} presented an approach for automated on-ramp merging and gap development considering vehicle speed constraints. Other efforts have also focused on multi-objective optimization problems for coordination of CAVs at intersection using either centralized or decentralized approaches \cite{Kamal2014, Makarem2013, qian2015}. 

Although previous research has aimed at enhancing our understanding of improving the efficiency through coordination of CAVs, deriving an optimal solution for a corridor  still remains a challenging control problem. 
In this paper, we address the problem of optimally coordinating CAVs that travel through a corridor under hard safety constraints, eliminating stop-and-go driving behavior. In previous work \cite{zhao2018, Zhao2019CCTA-1}, we presented a preliminary analysis on coordinating CAVs in a corridor without considering state, control, and safety constraints. In this paper, we implement a closed-form analytical solution that includes the rear-end safety constraint. 

The structure of the paper is organized as follows. In Section \ref{sec:pf}, we formulate the problem, and derive the analytical, closed form solution for the corridor with interior constraints. In Section \ref{sec:sim}, we validate the effectiveness of the analytical solution in a simulation environment and conduct a comparison analysis with traditional human-driven vehicles. Finally, we provide concluding remarks in Section \ref{sec:conc}.

\section{Problem Formulation} \label{sec:pf}
\begin{figure}[!ht]
\centering
\includegraphics[scale = 0.27]{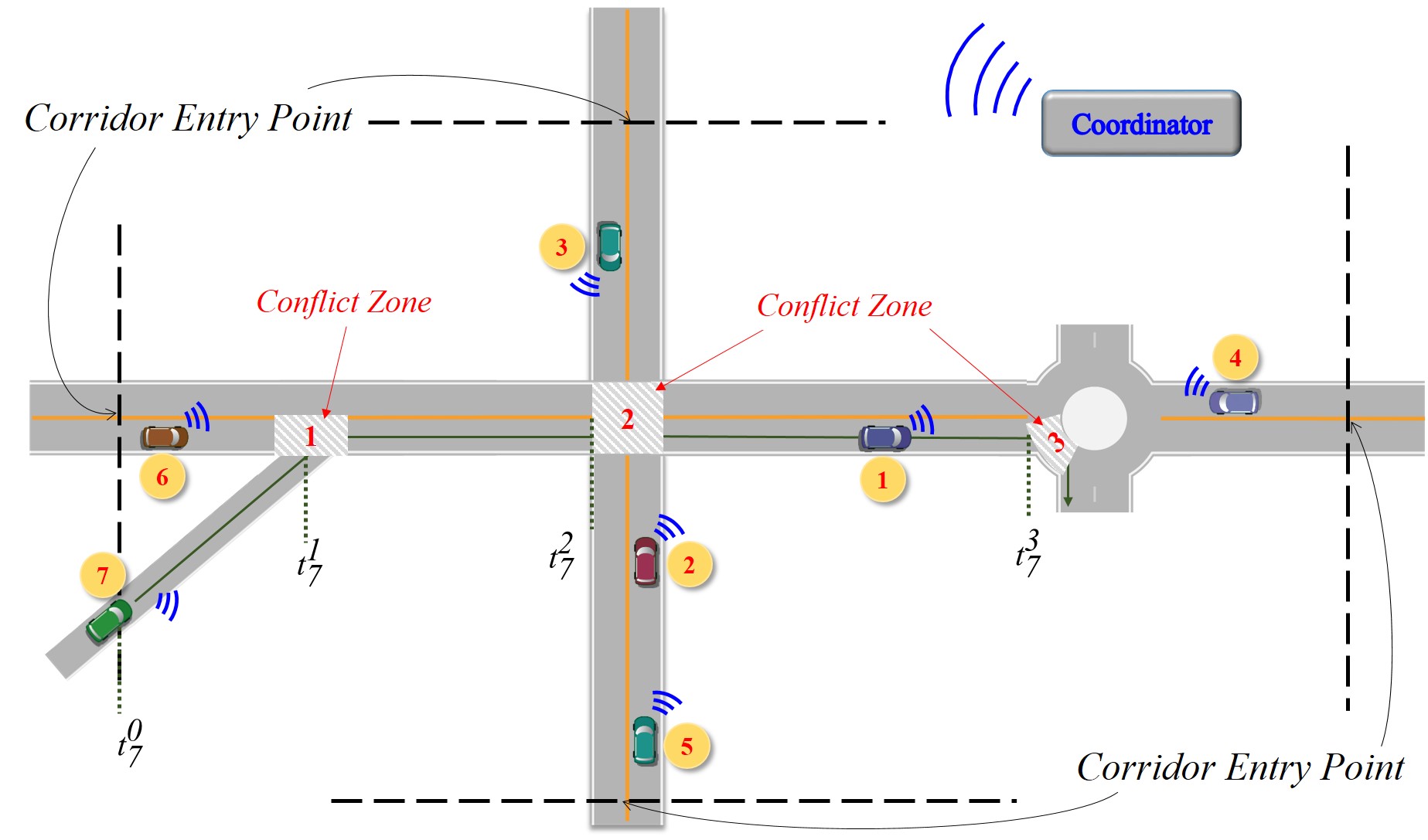} \caption{Corridor with
connected and automated vehicles.}%
\label{fig:corridor}%
\end{figure}

We consider a corridor (Fig. \ref{fig:corridor}) that consists of several conflict zones (e.g., a merging area, an intersection, and a roundabout), where potential lateral collision of vehicles may occur. The corridor has a coordinator that can monitor CAVs traveling along the corridor. Note that the coordinator is not involved in any decision on the CAV operation. 
The communication range of the coordinator can be adjustable and its length could be extended as needed. For example, we could use a network of drones to act as coordinators and broadcast with the CAVs. 

Let $N(t)\in\mathbb{N}$ be the number of CAVs in the corridor at time $t\in\mathbb{R}^{+}$, $\mathcal{N}(t)=\{1,2,\dots,N(t)\}$ be a queue of CAVs inside the corridor, and $\mathcal{Z}\in\mathbb{N}$ be the number of conflict zones along the corridor where lateral collisions may occur. When a CAV enters the boundary of the corridor, it broadcasts its route information to the coordinator. Then, the coordinator assigns a unique integer $i\in\mathbb{N}$ that serves as identification of CAVs inside the corridor. Let $t_i^0$ be the initial time that vehicle $i$ enters the corridor, $t_{i}^{z}$ be the time for vehicle $i$ to enter the conflict zone $z$, $z \in \mathcal{Z}$, and $t_{i}^{f}$ be the time for vehicle $i$ to enter the final conflict zone. For example, for CAV $\# 7$ (Fig.~\ref{fig:corridor}), $t_{7}^{0}$ is the time that it enters the corridor, $t_{7}^{1}$, $t_{7}^{2}$, and $t_{7}^{f}=t_{7}^{3}$ are the times that it enters the conflict zones $\#1, \#2,$ and $\#3$ respectively. There is a number of ways to assign $t_{i}^{z}$ for each CAV $i$. The policy through which the \textquotedblleft schedule\textquotedblright~  is specified is the result of a high-level optimization problem \cite{Malikopoulos2019CDC}.

The policy, which  determines the time $t_{i}^{z}$ that each CAV $i$ enter the conflict zone, can aim at maximizing the throughput at the corridor while ensuring that any lateral collision constraint never becomes active. On the other hand, for each CAV $i$, deriving the optimal control input (minimum acceleration/deceleration) to achieve the target $t_{i}^{z}$ can aim at minimizing its fuel consumption \cite{Malikopoulos2010a} while ensuring that the rear-end collision avoidance constraint never becomes active.
In what follows, we consider that a control scheme for determining
$t_{i}^{z}$ for each CAV $i$ is given, and we will focus on the low-level control problem that yields for each CAV the optimal control input to achieve the assigned $t_{i}^{z}$ subject to the state, control, and rear-end collision avoidance constraints.

\subsection{Vehicle model, Constraints, and Assumptions}
Each CAV $i\in\mathcal{N}(t)$ is modeled by a second order dynamics
\begin{equation}%
\dot{p}_{i} = v_{i}(t), \quad
\dot{v}_{i}= u_{i}(t),
\label{eq:model2}
\end{equation}
where $p_{i}(t)\in\mathcal{P}_{i}$, $v_{i}(t)\in\mathcal{V}_{i}$, and $u_{i}(t)\in\mathcal{U}_{i}$ denote the position, speed and control input (acceleration/deceleration) of each CAV $i$ in the corridor; $s_{i}(t)\in\mathcal{S}_{i}$ denotes the distance between CAV $i$ and CAV $k$ which is directly ahead of $i$, and $\xi_{i}$ is the reaction constant of the CAV. The sets $\mathcal{P}_{i}$, $\mathcal{S}_{i}$, $\mathcal{V}_{i}$, and $\mathcal{U}_{i}$,  $i\in\mathcal{N}(t),$ are complete and totally bounded subsets of $\mathbb{R}$. In the rest of the paper, we reserve the symbol $k$ to denote the CAV which is physically immediately ahead of $i$ in the same lane.

Let $x_{i}(t)=\left[p_{i}(t) ~ v_{i}(t)\right]  ^{T}$ denote the state of each CAV $i$, with initial value $x_{i}^{0}=\left[p_{i}^{0} ~ v_{i}^{0}\right]  ^{T}$, where $p_{i}^{0}= p_{i}(t_{i}^{0})=0$ and $v_{i}^{0}= v_{i}(t_{i}^{0})$. 
To ensure that the control input and speed are within a given admissible range, the following constraints are imposed.
\begin{equation}%
\begin{split}
u_{min} &  \leq u_{i}(t)\leq u_{max},\quad\text{and}\\
0 &  \leq v_{min}\leq v_{i}(t)\leq v_{max},\quad\forall t\in\lbrack t_{i}%
^{0},t_{i}^{f}],
\end{split}
\label{speed_accel_constraints}%
\end{equation}
where $u_{min}$, $u_{max}$ are the minimum deceleration and maximum
acceleration for each CAV $i\in\mathcal{N}(t)$, and $v_{min}$, $v_{max}$ are the minimum and maximum speed limits respectively. 
To ensure the absence of rear-end collision of CAV $i \in \mathcal{N}(t)$ and its immediately preceding CAV $k \in \mathcal{N}(t)$, we impose the following rear-end safety constraint 
\begin{equation}
\begin{split}
s_{i}(t)=\xi_i \cdot (p_{k}(t)-p_{i}(t)) \ge \delta_i(t),~ \forall t\in [t_i^0, t_i^f].
\label{eq:rearend}
\end{split}
\end{equation}
Here, the minimum safe distance $\delta_i(t)$ is a function of speed $v_i(t)$ since $\delta_i(t)=\gamma_i + \rho_i \cdot v_i(t)$, for all $t\in [t_i^0, t_i^f]$, 
where $\gamma_i$ is the standstill distance, and $\rho_i$ is minimum time gap that CAV $i$ would maintain while following another CAV. 
In the  modeling framework described above, we impose the following assumptions:
\begin{assumption} \label{ass1}
All CAVs are connected and automated. 
\end{assumption}

\begin{assumption} \label{ass2}
For each CAV $i$, none of the constraints \eqref{speed_accel_constraints} and \eqref{eq:rearend} is active at $t^0_i$.
\end{assumption}

\begin{assumption}\label{ass3}
Each CAV $i$ has proximity sensors and can measure local information without errors or delays.
\end{assumption}

\begin{assumption}\label{ass4}
The communication between the coordinator and each CAV $i$ is reliable and instantaneous.
\end{assumption}

\begin{assumption} \label{ass: no_lane_change}
The corridor only contains single-lane road segments. The CAVs traveling in the corridor do not change lanes except to make necessary turns.
\end{assumption}

The first assumption limits the scope of our paper to an idealized environment where all vehicles are connected and automated, e.g., 100\% CAV penetration rate.
The second assumption ensures that the solution of the optimal control problem starts from a feasible state and control input. 
The third assumption might impose barriers during implementation, but can be relaxed if the noise in the measurements and delays are bounded.
The fourth assumption is necessary for the upper-level optimization problem although this problem is not addressed in this paper. 
The last assumption simplifies the optimal control problem so as to avoid implications related to lane changing. However, the proposed framework could be extended to multiple lanes by enhancing the vehicle model to account for this accordingly. 

\section{The Low-Level Optimal Control Problem} \label{sec:sol}

We consider the problem of minimizing the control input for each CAV $i\in\mathcal{N}(t)$ from the time $t_i^0$
that the CAV $i$ enters the control zone until the time $t_i^{f}$ that it exits the last conflict zone under the hard safety constraint to avoid rear-end collision. By minimizing each CAV's control input, we minimize transient engine operation. Thus, we can have direct benefits in fuel consumption and emissions since internal combustion engines are optimized over steady state operating points (constant torque and speed) \cite{Malikopoulos2008sae}.

Therefore, the optimization problem for each CAV $i\in\mathcal{N}(t)$ is formulated as follows:
\begin{gather}\label{eq:decentral}
J_{i}(u(t))=  \frac{1}{2} \int_{t^0_i}^{t^f_i} u^2_i(t)~dt,\\ 
\text{subject to:}
~\eqref{eq:model2},\eqref{eq:rearend},\nonumber\\
\text{and given }t_{i}^{0}\text{, }v_{i}^{0}\text{, }\text{
}p_{i}(t_{i}^{0})=0\text{, } {t_i^{z}}^*, p_{i}({t_i^{z}}^*)=p_z.\nonumber
\end{gather}
Note that we have omitted the state and control constraint \eqref{speed_accel_constraints}. The  problem formulation with the state 
and control constraints requires the constrained and unconstrained arcs of the state and control input to be pieced together to satisfy the Euler-Lagrange equations and necessary condition of optimality. To simplify the analysis, we focus on the rear-end safety constrained case \eqref{eq:rearend} only. The other constrained cases related to the state and control are similar to the cases presented in \cite{Malikopoulos2017,Malikopoulos2019ACC,Mahbub2020ACC-1}, and thus, we do not repeat here.
In our analysis, we consider that when a CAV enters the control zone, none of the constraints is active (Assumption \ref{ass2}). From \eqref{eq:model2} and \eqref{eq:decentral}, we formulate the Hamiltonian for each CAV $i\in\mathcal{N}(t)$ by indirectly adjoining the rear-end safety constraint \eqref{eq:rearend},
\begin{gather}
H_{i}\big(t, p_{i}(t), v_{i}(t), s_{i}(t), u_{i}(t)\big)  \nonumber \\
=\frac{1}{2} u(t)^{2}_{i} + \lambda^{p}_{i} \cdot v_{i}(t) +\lambda^{s}_{i} \cdot \xi_i \cdot (v_{k}(t) - v_{i}(t)) \nonumber + \lambda^{v}_{i} \cdot u_{i}(t) \\
+ \mu^{e}_{i} \cdot (\rho_i \cdot u_i(t)-\xi_i(v_k(t)-v_i(t))) ,
\label{eq:hamil}
\end{gather}
where $\lambda^{p}_{i}$, $\lambda^{s}_{i}$, and  $\lambda^{v}_{i}$ are the costate components, and $\mu^{e}_i$ is the Lagrange multiplier.
The Euler-Lagrange equations become
\begin{gather}
\dot\lambda^{p}_{i}(t) = - \frac{\partial H_i}{\partial p_{i}} = 0, \quad \dot\lambda^{s}_{i}(t) = - \frac{\partial H_i}{\partial s_{i}} = 0,\nonumber \\
\dot\lambda^{v}_{i}(t) = - \frac{\partial H_i}{\partial v_{i}} = -(\lambda^{p}_{i} - \xi_i \cdot \lambda^{s}_{i} + \mu_i^e\cdot \xi_i),~ \text{and} \nonumber \\
\frac{\partial H_i}{\partial u_{i}} = u_{i}(t) + \lambda^{v}_{i} + \rho_i\cdot \mu_i^e = 0. \label{eq:ELp}
\end{gather}

\subsubsection{State constraint is not active, analytical solution without interior constraints} \label{case1}
When the rear-end constraint is not active, $\mu^{e}_{i}=0$. 
From (\ref{eq:ELp}), we have $\lambda^{p}_{i}(t) = a_{i}$,   
$\lambda^{s}_{i}(t)= b_{i}$, and $\lambda^{v}_{i}(t) = -\big((a_{i} - b_{i} \cdot \xi_i)\cdot t + c_{i}\big)$. 
The coefficients $a_{i}$, $b_{i}$, and $c_{i}$
are constants of integration corresponding to each CAV $i$. From \eqref{eq:ELp}, the optimal control input as a function of time is
\begin{equation}
u^{*}_{i}(t) = (a_{i} - b_{i} \cdot \xi_i) \cdot t + c_{i}, ~ \forall t \ge t^{0}_{i}. \label{eq:24}
\end{equation}

Substituting the last equation into \eqref{eq:model2} we find the optimal speed and position for each CAV,
namely
\begin{gather}
v^{*}_{i}(t) = \frac{1}{2} (a_{i} - b_{i} \cdot \xi_i) \cdot t^2 + c_{i} \cdot t +d_{i}, ~ \forall t \ge t^{0}_{i}\label{eq:25}\\
p^{*}_{i}(t) = \frac{1}{6} (a_{i} - b_{i} \cdot \xi_i) \cdot t^3 +\frac{1}{2} c_{i} \cdot t^2 + d_{i}\cdot t +e_{i}, ~ \forall t \ge t^{0}_{i}. \label{eq:26}%
\end{gather}
where $d_{i}$ and $e_{i}$ are constants of integration. The constants of integration $a_i$, $c_{i}$, $d_{i}$, and $e_{i}$ are computed for each CAV $i$ from the values of the speed and position at $t_i^0$, the position at $t_i^f$, and $\lambda^{v}_{i}(t^{f}_{i})=0$.

\subsubsection{State constraint is active, analytical solution without interior constraints}\label{case2}
Suppose the CAV starts from a feasible state and control at $t=t_{i}^0$ and at some time $t=t_{1}$, $s_{i}(t_{1})=\delta(t_{1})$. In this case, $\mu_{i}^{e} \neq 0$.
The closed-form solution for this case has been derived in \cite{Malikopoulos2019ACC}, and thus, it is not repeated here.

\subsubsection{State constraint is not active, analytical solution with interior constraints}\label{case3}
In this case, the path of CAV $i$ consists of more than one conflict zone, e.g., CAV $i$ enters from the second ramp and travels through \#1 and \#2 in Fig. \ref{fig:corridor}. Between the time $t_i^0$ that the CAV enters the corridor and the time $t_i^f$ that the CAV exits the conflict zone \#2, CAV $i$ has to travel across the intermediate merging zone \#1 (i.e., the first conflict zone) at designated time ${t_i^{1}}^*$, then we have an additional interior boundary condition at $t = {t_i^{1}}^*$, $p_i({t_i^{1}}^*) = p_1$ and/or $v_i({t_i^{1}}^*) = v_1$, where where $p_1$ and $v_1$ are the position and designed speed limit (if exists) for the conflict zone. In this case, the Hamiltonian in \eqref{eq:hamil} needs to be solved as a three point boundary value problem piecing two unconstrained arcs together with proper jump conditions at the interior boundary point at $t = {t_i^{1}}^*$ as described in \cite{Mahbub2019ACC}.
Following the same procedure as in \cite{Mahbub2019ACC}, we can extend this solution for multiple interior points of a route of a CAV.

\subsubsection{State constraint is active, analytical solution with interior constraints} \label{case4}
We consider the case where the path of CAV $i$ consists of more than one conflict zone (i.e., analytical solution with interior constraints), and the safety constraint is activated in one or more  arcs. In this case, we have at least two junction points where the Hamiltonian is discontinuous. Suppose the CAV starts from a feasible state and control at $t=t_{i}^0$, travel across an intermediate conflict zone $1$ at a given time $t={t_i^{1}}^*$, and at some time $t=t_{1}$, $s_{i}(t_{1})=\delta(t_{1})$. Note that safety constraint may be violated before or after ${t_i^{1}}^*$.
Depending on relationship between ${t_i^{1}}^*$ and $t_{1}$, the constrained and unconstrained arcs can be pieced together at the optimal junction points by combining case 2 and 3.

\section{Simulation Results} \label{sec:sim}
\subsection{Numerical Analysis}
We analyzed four cases to validate the effectiveness of the analytical solution. In all cases, we set a leading CAV and the following CAV, we call it ``ego" vehicle, traveling through a corridor. The total length of the corridor is 300 m. When $t_i^0=0$ $s$, the following CAV $i$ is at the entry of the corridor with $p_i^0=0$ $m$. For simplification, we assume the final time for the following CAV has been determined at $t_i^f = 26$ $s$. At $t_i^0=0$ $s$, the speed for the leading CAV is $v_k=11.5$ $m/s$. The final time of the leading CAV around $t_f^k=24$ $s$, after which, the leading CAV exits the corridor. We use different acceleration profiles of the leading CAVs to validate the optimal control for the following CAV.

\begin{figure} [ht]
	\centering
	\includegraphics[width=0.45\textwidth]{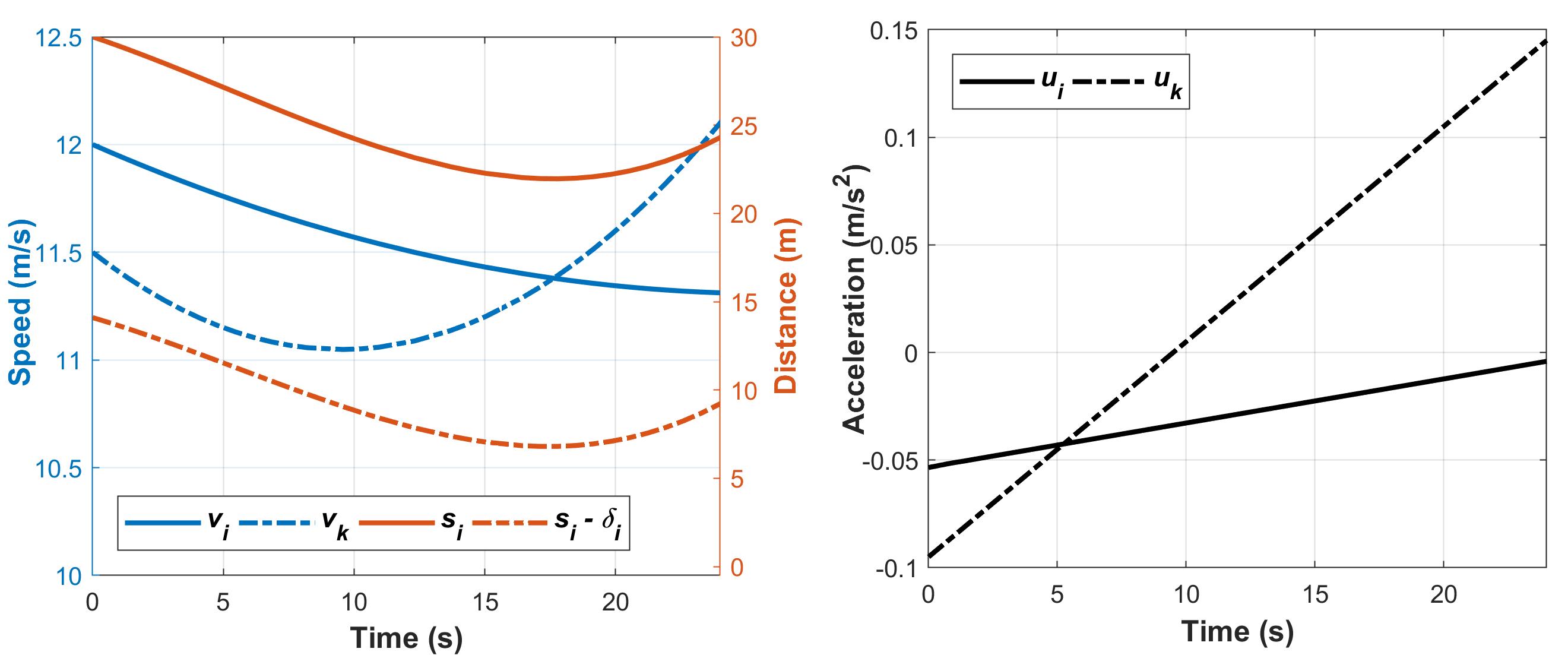}   
	\caption{Case 1: safety constraint is not activated, no interior points.} 
	\label{fig:case1a}
\end{figure}

Case 1: safety constraint not activated, no interior points. 
In this case, CAV $i$ enters with an initial speed of 12.0 $m/s$ at time $t_i^0=0$ $s$. The initial following distance $s_i^0$ is 30 $m$. The state following distance $s_i(t)$, together with the difference between following distance $s_i(t)$ and the safety distance $\delta_i(t)$, is plotted in the left panel of Fig. \ref{fig:case1a}, where we observe that if state constraint is not activated, the optimal control yields a linear acceleration profile for the following CAV $i$. 

\begin{figure} [ht]
	\centering
	\includegraphics[width=0.45\textwidth]{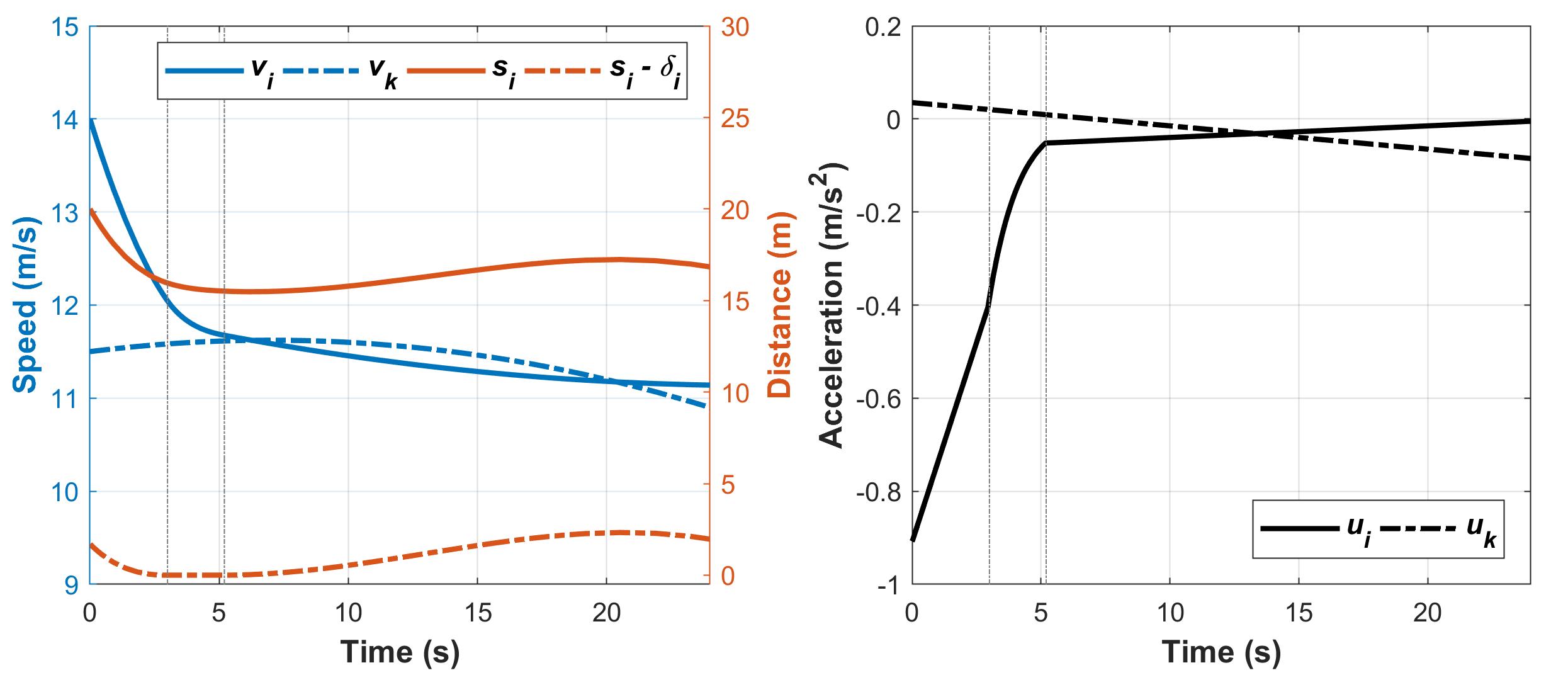}   
	\caption{Case 2: safety constraint is activated, no interior points.} 
	\label{fig:case1b}
\end{figure}

Case 2: safety constraint activated, no interior points.
In this case, CAV $i$ enters with a higher initial speed of 14.0 $m/s$ at time $t_i^0=0$ $s$, when the initial following distance $s_i^0$ is 20 $m$. We see from Fig. \ref{fig:case1b} that the state constraint becomes activated at $t_1=3.2$ $s$ when $s_i(t_1) - \delta_i(t_1)=0$. Since CAV $k$ keeps accelerating while CAV $i$ decelerates after $t_1$, CAV $i$ exits the constrained arc at $t_2=5.2$ $s$. After $t_2$, CAV $i$ decelerates with a linearly increased acceleration profile.

\begin{figure} [ht]
	\centering
	\includegraphics[width=0.45\textwidth]{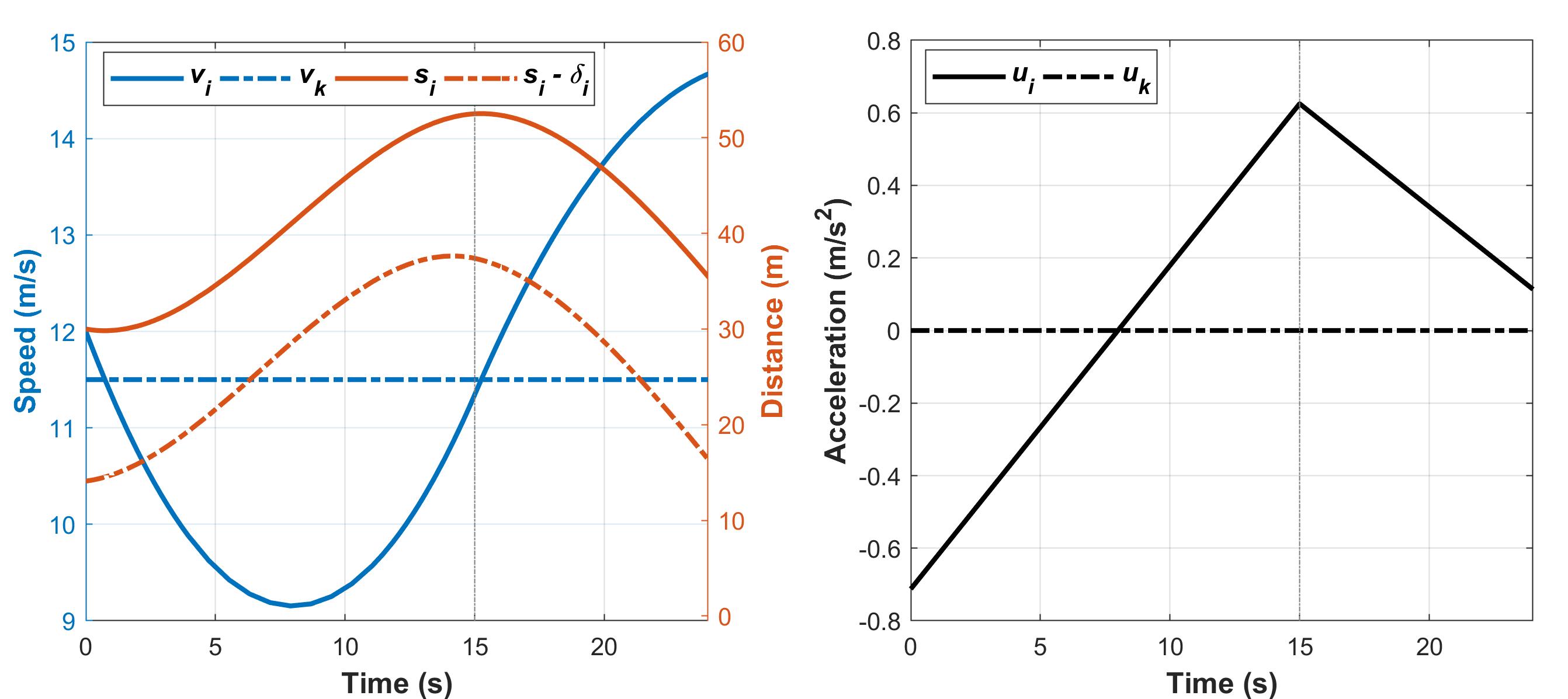}   
	\caption{Case 3: safety constraint is not activated, interior point exists.} 
	\label{fig:case1c}
\end{figure}

Case 3: safety constraint not activated, interior point exists.
We look into the case when there is an interior point, the case is set up as follows: CAV $i$ enters with an initial speed of 12.0 $m/s$ at time $t_i^0=0$ $s$, when the initial following distance $s_i^0=0$ $m$. There is an intersection at $p_1=150$ $m$, where the scheduled entry time for CAV $i$ is $t_1=15$ $s$. Suppose CAV $k$ travels at a constant speed. Due to the existence of an intermediate point, we can see two arcs from the optimal acceleration profile the right panel of Fig. \ref{fig:case1c}. For the first half segment, CAV $i$ first decelerate and then accelerate to meet the assigned entry time at the intersection; then for the second half, CAV $i$ keeps accelerating till the end.

\begin{figure} [ht]
	\centering
	\includegraphics[width=0.45\textwidth]{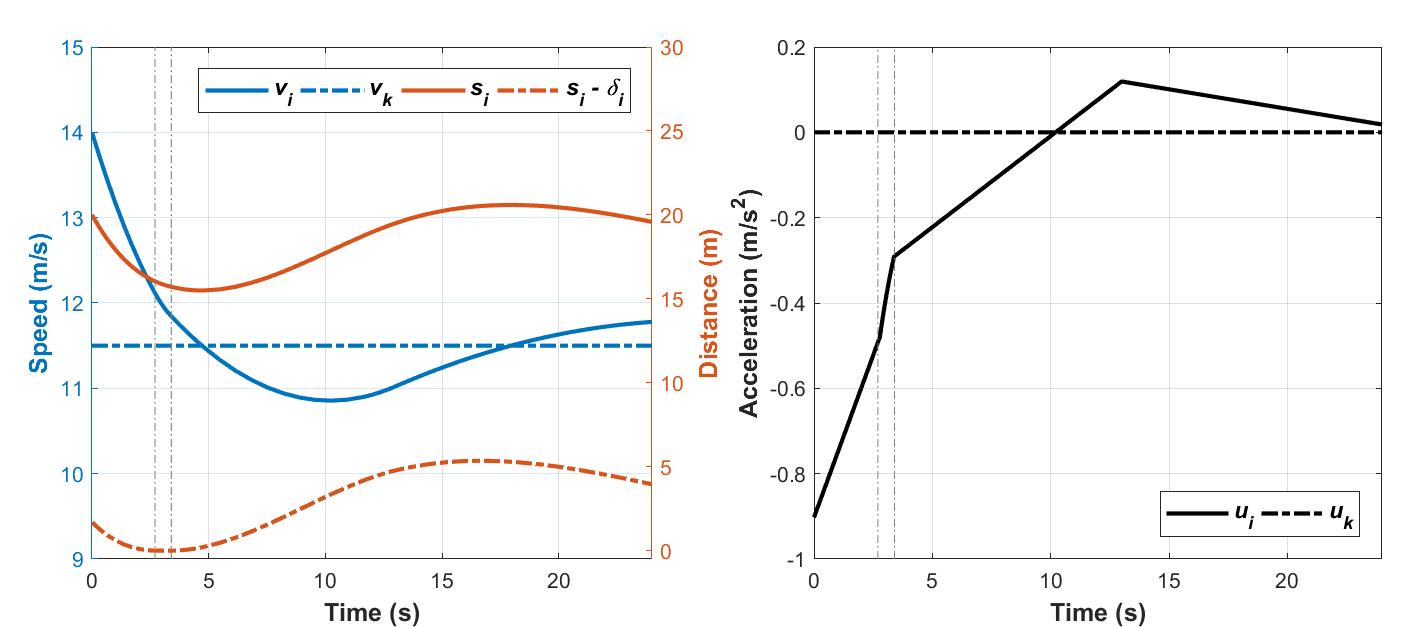}   
	\caption{Case 4: safety constraint is activated, interior point exists.} 
	\label{fig:case1d}
\end{figure}

Case 4: safety constraint is activated, interior point exists.
In case 3, if the assigned entry time of CAV $i$ at the intermediate point, the safety constraint may be activated. Thus, in case 4, we analyze the situation when safety constraint is activated in a corridor. In this case, the scheduled entry time at $p_1=150$ $m$ for CAV $i$ is $t_1=13$ $s$. In Fig. \ref{fig:case1d}, we see that due to the change of entry time at the intermediate point, CAV $i$, with a high initial speed, activates the safety constraint at 2.7 $s$, and leaves the constrained arc at 3.4 $s$. With the optimal control for constrained and unconstrained arcs, CAV $i$ is able to pass through the intersection and final location at predetermined times.

\subsection{Traffic Simulation}
\begin{figure} [!ht]
	\centering
	\includegraphics[scale=0.56]{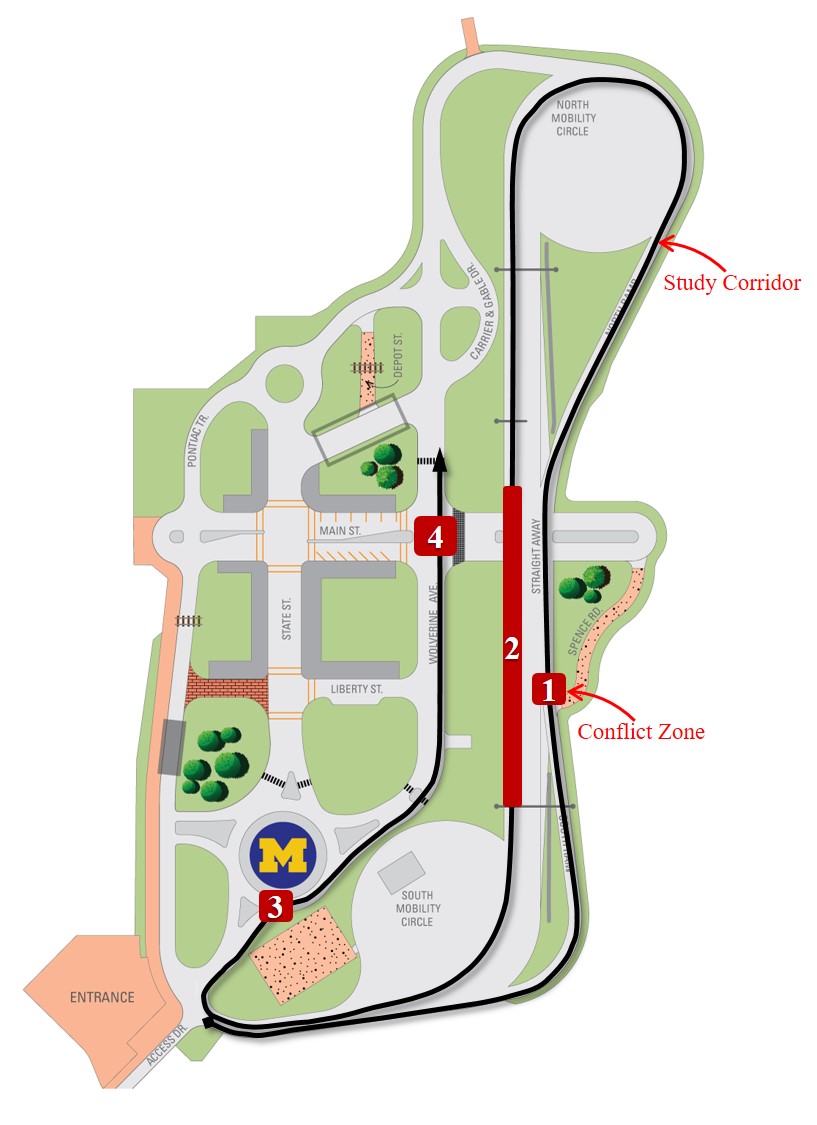}   
	\caption{Study corridor in Mcity.} 
	\label{fig:mcity}
\end{figure}

With the simulation network of Mcity created in PTV VISSIM (Version 11) environment, we define a corridor consisting of four conflict zones: (1) a merging roadway, (2) a speed reduction zone, and (3) a roundabout, and (4) an intersection (Fig. \ref{fig:mcity}). CAVs enter the network on the ramp, join the traffic on the highway with desired speed of 22 $m/s$, and then enter the speed reduction zone where the speed limit drops to 11 $m/s$. The CAVs exit the highway segment and travel through the roundabout, where a desired speed of 13 $m/s$ is imposed until the exit of the roundabout, to the intersection (conflict zone $\#4$). 

To evaluate the network performance with the proposed control framework, we consider two scenarios as follows: \textbf{(a) Baseline:} All vehicles in the network are non-connected and non-automaicles. In this case, the Wiedemann car following model \cite{wiedemann1974} built in VISSIM is applied. 1.2 $s$ time headway is adopted to estimate the minimum allowable following distance.\newline
\textbf{(b) Optimal Control:} We adopt the same simulation platform as in our previous work \cite{Malikopoulos2018a}, where the proposed control framework is integrated to generate the optimal acceleration/deceleration profile for each CAV in the network. Same time headway under Scenario 1 is applied in the optimal control model.

\subsubsection{Vehicle Trajectory}
\begin{figure} [!ht]
	\centering
	\includegraphics[width=0.5\textwidth]{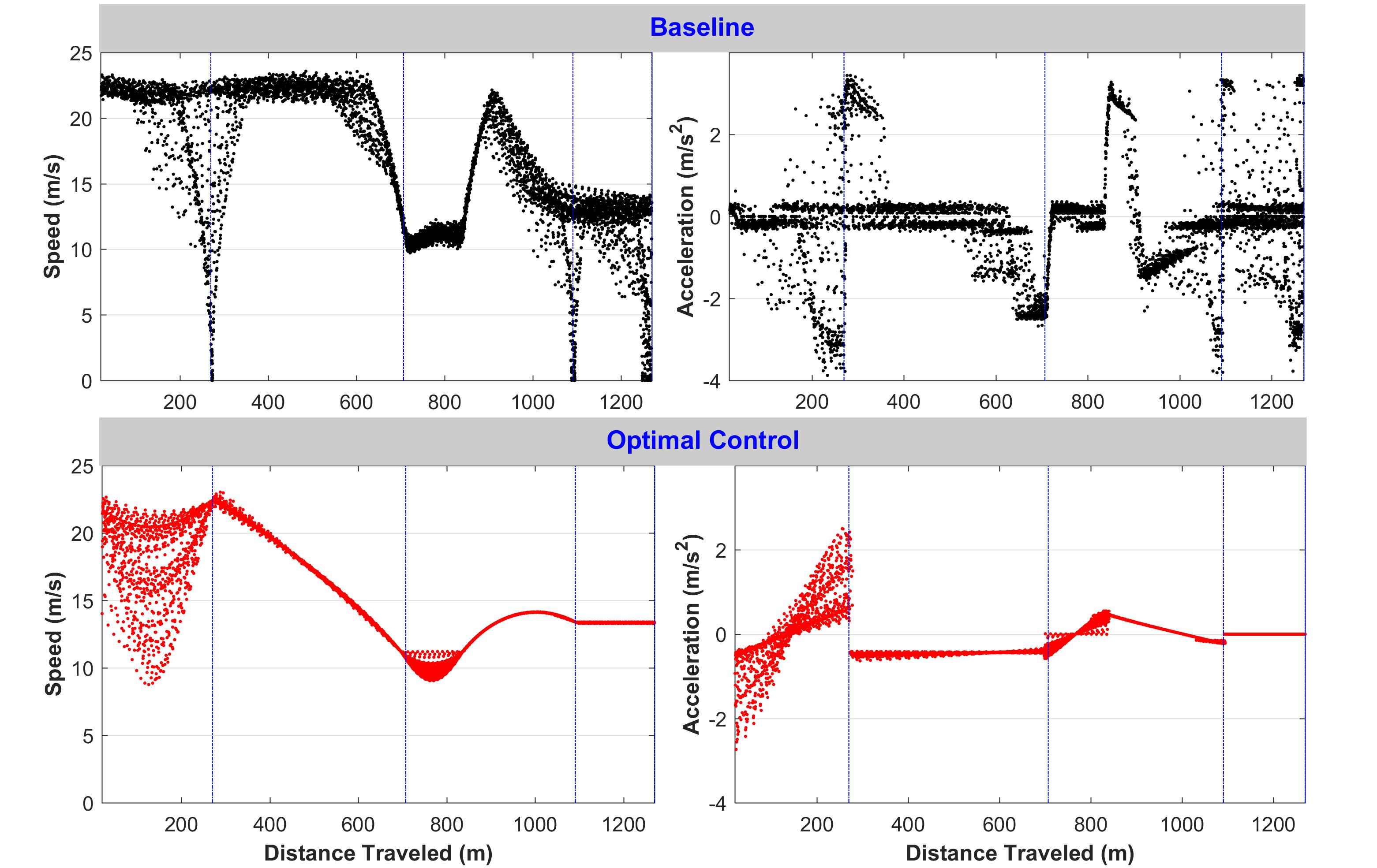}   
	\caption{Vehicle trajectories.} 
	\label{fig:trajectory}
\end{figure}

The speed trajectories for 0\% and 100\% CAV penetration rate are shown in Fig. \ref{fig:trajectory}. For the baseline scenario with 0\% CAV penetration rate, CAVs traveling along the corridor need to yield to mainline traffic, and wait for the green light before the intersection. Thus, we can see that there are many fluctuations in CAV acceleration profiles under the baseline scenario (i.e., black dots in upper right panel of Fig. \ref{fig:trajectory}). With the optimal control operation under 100\% CAV penetration rate, the traffic information for the entire corridor is shared for all CAVs. Therefore, CAVs traveling through the corridor could drive more smoothly to avoid hard acceleration/deceleration for any merging or speed reduction events in the path (i.e., red dots in lower panel of Fig. \ref{fig:trajectory}). Note that a desired speed is defined for each conflict zone, which is the reason for the jumps in the acceleration profiles. If desired, we could eliminate the jumps in the optimal acceleration solution by removing the hard speed limits for the conflict zones (as shown in Fig. \ref{fig:case1c}).

\subsubsection{Travel Time}
\begin{figure} [ht]
	\centering
	\includegraphics[width=0.5\textwidth]{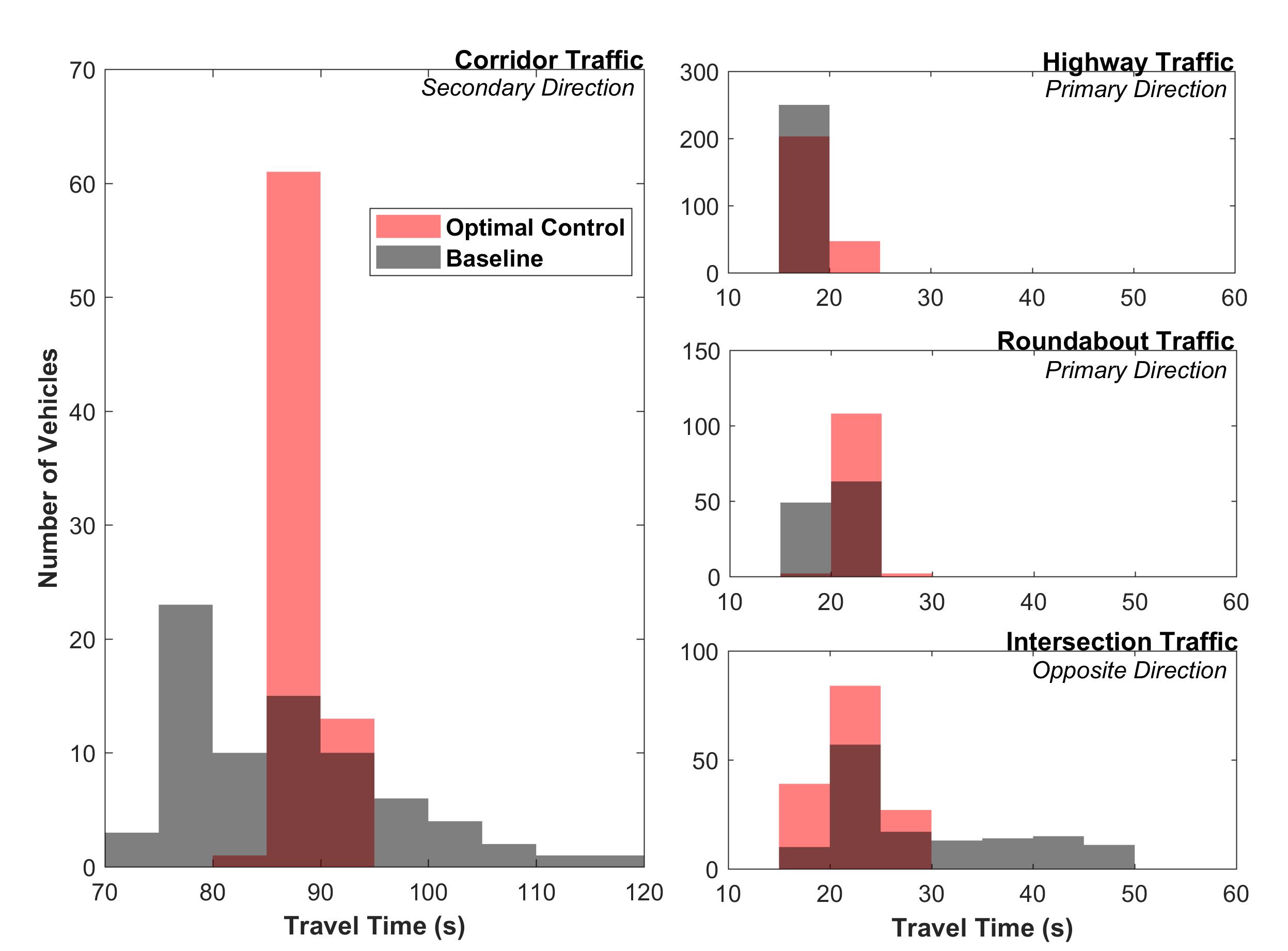}   
	\caption{Travel time distribution for different network segments.} 
	\label{fig:travel_time}
\end{figure}

As we can see from the left panel of Fig. \ref{fig:travel_time}, smooth vehicle movement under the optimal control operation leads to smaller fluctuation in corridor travel time among CAVs, compared with baseline scenario. The right three panels in Fig. \ref{fig:travel_time} represent travel time distribution for mainline traffic towards the three conflict zones, i.e., highway merging area, roundabout, and the intersection. In general, the average travel time for mainline CAVs is scarified in order to generate proper gaps for CAVs on secondary direction to conduct non-stop merging maneuver. However, since traffic signal controller is disabled in the network with 100\% CAV penetration rate, the travel times for all CAVs towards the intersection are reduced substantially (i.e., the bottom figure in the right panel of Fig. \ref{fig:travel_time}).

\subsubsection{Following Distance}
\begin{figure} [ht]
	\centering
	\includegraphics[width=0.5\textwidth]{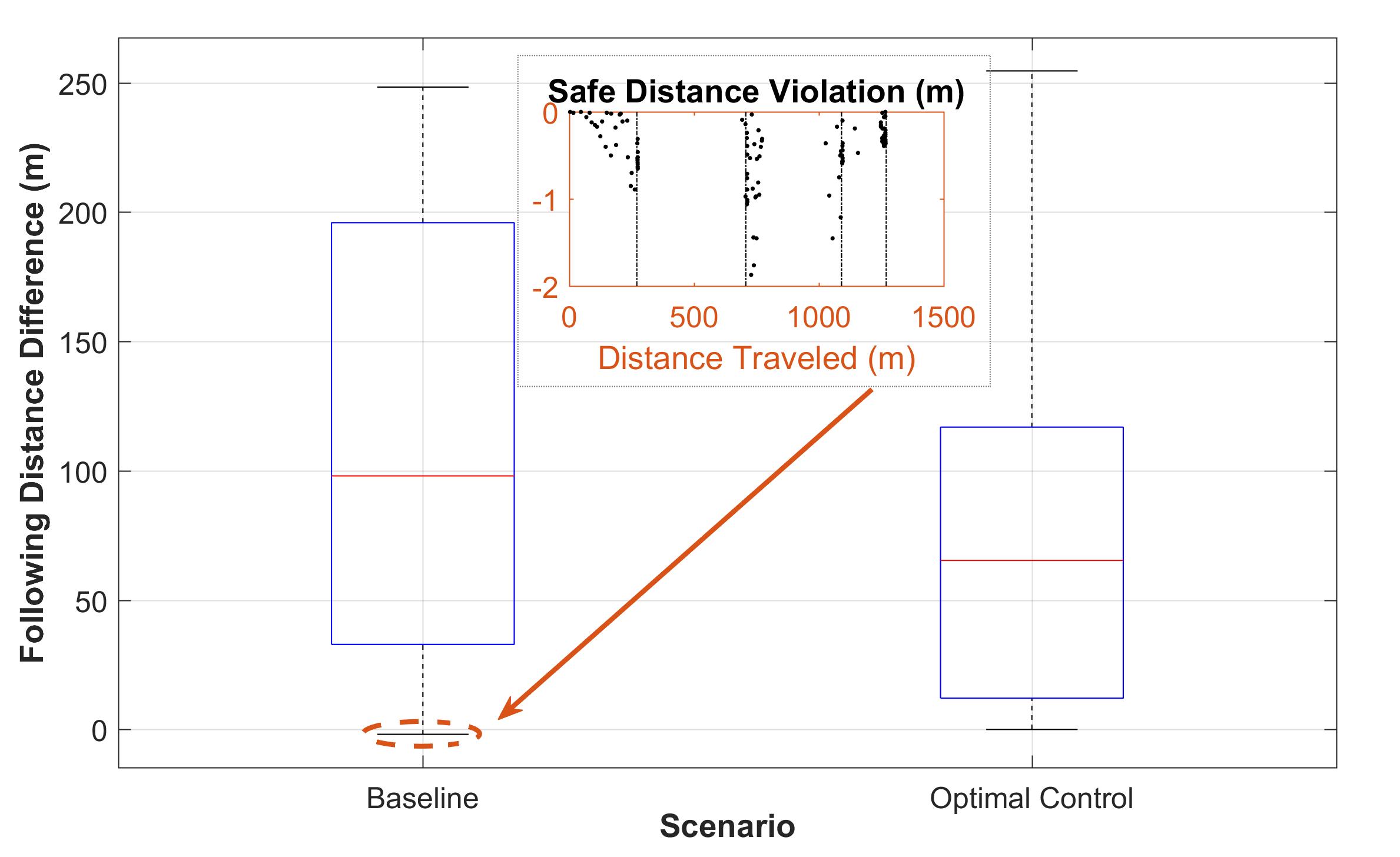}   
	\caption{Accumulated delay in each zone.} 
	\label{fig:safe_distance}
\end{figure}

We collect from VISSIM the distance between a CAV and its physically immediately leading CAV (i.e., following distance) every second, and the minimum safe distance calculated based on travel speed. The following distance difference in Fig. \ref{fig:safe_distance} is defined as the difference between second-by-second CAV following distance and the minimum safe distance. A value above 250 $m$ in VISSIM records means that there is no leading CAV ahead. 
We see in Fig. \ref{fig:safe_distance} that with the optimal control algorithm, the variance in following distance difference among CAVs decreases, indicating a more homogeneous traffic pattern under 100\%  penetration rate. Furthermore, we note that in the baseline scenario, there are some cases when safety constraint is violated as circled by the eclipse in Fig. \ref{fig:safe_distance}. Plotting the events along the distance traveled in the corridor (i.e., the insert panel in Fig. \ref{fig:safe_distance}), we see that most of the violation events happen near the entries of conflict zones as well as the start of the speed reduction zone. If sudden/sharp deceleration (e.g., from a relative high speed to a low speed or even stopping) is necessary for a series of CAVs, chances are that some CAVs in the chain may not be able to decelerate enough in a short time period with the minimum safe distance constraint satisfied (i.e., relatively high risk of collision). It implies that by recommending acceleration/deceleration profiles for all the CAVs, the optimal control algorithm could potentially reduce the risk of collision and improve traffic safety.

\subsubsection{Fuel Consumption}
With vehicle trajectory data collected every 1 $s$, fuel consumption is estimated by using the polynomial metamodel proposed in \cite{Kamal2011} that relates vehicle fuel consumption as a function of speed $v(t)$ and acceleration $u(t)$. 
With smooth acceleration/deceleration profiles throughout the entire corridor (see Fig. \ref{fig:trajectory}), vehicles' stop-and-go driving behavior is eliminated under Scenario 2 with 100\% CAV penetration rate. Thus, transient engine operation is minimized, leading to direct fuel consumption savings compared to the baseline scenario.
Overall, through the optimal control algorithm, an average of 41\% savings in total fuel consumption for CAVs traveling along the corridor is yielded.

\section{Concluding Remarks} \label{sec:conc}
In this paper, we investigated the optimal coordination of CAVs in a corridor.  We derived a closed-form analytical solution that considers interior constraints and provides the optimal trajectory for the entire route for each CAV. We showed through simulation that  coordination of CAVs can eliminate stop-and-go driving and improve fuel consumption. Ongoing work addresses the problem of optimal control of CAVs in a mixed traffic environment. 

Although potential benefits of full penetration rates of CAVs to alleviate traffic congestion and reduce fuel consumption are apparent, different penetration rates of CAVs can alter significantly the efficiency of the entire system. Ongoing work emphasizes addressing the problem of coordinating CAVs in a mixed traffic environment.
In our proposed framework, we made the assumption of perfect communication which might impose barriers in a potential implementation and deployment of the proposed framework. 
Another important direction for future research is to relax this assumption, and investigate the implications of having information with errors and/or delays to the system behavior.


\bibliographystyle{IEEEtran}
\bibliography{corridor_ref}

\end{document}